\voffset=1cm

\input amssym.def
\input amssym.tex

\hsize=125mm
\vsize=195mm

\input psfig

\parindent=0pt
\magnification=1100
\baselineskip=18pt
\nopagenumbers

\input latexpicobjs
\input pictex

\font\titel=cmbx10 at 20 truept
 at 20 truept
 at 20 truept
 at 14 truept
\font\bfmittel=cmbx8 at 14 truept
 at 12 truept
\font\klein=cmr8

\font\bfklein=cmbx8

 at 10 truept

\font\refklein=cmr12 at 11 truept
\font\titel=cmbx10 scaled\magstep2

\def\loba{\hbox{J\kern-1pt I}}

\def\gd{\hbox{$\bullet\kern-2pt$ --\kern-4pt---\kern-4pt--- 
\kern-2pt}}
\def\graphd{\hbox{$\bullet\kern-2pt$ --\kern-4pt---\kern-4pt---
$\kern-2pt\circ$}}
\def\gr #1{\hbox{$\bullet\kern-2pt $ ---\kern-3pt\raise
5pt\hbox{$#1$}\kern-7pt --\kern-6pt ---\kern-5pt ---\kern+1pt}}
\def\grh #1{\hbox{$\bullet\kern-2pt $ ---\kern-3pt\raise
6pt\hbox{$#1$}\kern-7pt --\kern-6pt ---\kern-5pt ---\kern+1pt}}

\def\r#1{\hbox{$#1$}\raise-8pt\hbox{${\,\tilde{}}$}}

\def\Vor{\parindent=17pt\par\hang\textindent}
\def\ad#1{#1}

\nopagenumbers
\def\makeheadline{\vbox to 0pt{\vskip-37pt\line{\vbox to
8.5pt{}\the\headline}\vss}\nointerlineskip}
\headline={\ifodd\pageno\rightheadline\else\leftheadline\fi}
\def\makeheadline{\vbox to 0pt{\vskip-37pt\line{\vbox to
8.5pt{}\the\headline}\vss}\nointerlineskip}
\headline={\ifnum\pageno=1\hfil\else\kopf\fi}
\def\kopf{\ifodd\pageno\rightheadline\else\leftheadline\fi}
\def\rightheadline{\tenrm\hfil On the growth of cocompact hyperbolic Coxeter groups \hfil\folio}
\def\leftheadline{\tenrm\folio\hfil Ruth Kellerhals and Genevi\`eve Perren\hfil} 
\voffset=2\baselineskip


\centerline{{\titel On the growth of}}
\medskip
\centerline{{\titel cocompact hyperbolic Coxeter groups}}
\vskip1cm
{\hfill{Ruth KELLERHALS*\quad and\quad Genevi\` eve PERREN}\footnote*{Partially 
supported by Schweizerischer Nationalfonds 200020-121506/1, 200020-113199/1.

}\hfill}
\smallskip{\hfill{\klein ruth.kellerhals@unifr.ch\hskip1.6cm genevieve.perren@unifr.ch}\hfill}

{\hfill{\klein Department of Mathematics, University of Fribourg, Switzerland}\hfill}
$$ $$
\centerline{
\vbox{
\hsize=130 true mm
\noindent {\baselineskip=14pt
{\bfklein Abstract.}\quad 
\klein 
For an arbitrary cocompact hyperbolic Coxeter group $G$
with finite generator set $S$ and complete growth function $\,f_S(x)=P(x)/Q(x)\,$, we provide a recursion
formula for the coefficients of the denominator polynomial $Q(x)$. It allows to 
determine recursively the Taylor coefficients and to study
the arithmetic nature of the poles of the growth function $f_S(x)$ in terms
of its subgroups and exponent variety. We illustrate this in the easy case of compact right-angled hyperbolic $n$-polytopes.
Finally, we provide detailed insight into the case of
Coxeter groups with at most 6 generators, acting cocompactly on hyperbolic 4-space,
by considering the three combinatorially different families
discovered and classified by Lann\'er, Kaplinskaya and Esselmann, respectively. 
}}}
$$ $$

{\bfmittel 1. Overview and results}
\bigskip
\medskip
Let $G$ be a discrete group generated by finitely many reflections in hyperplanes (mirrors) of hyperbolic space $\Bbb H^ n$
such that the orbifold $\Bbb H^ n/G\,$  is compact. We call $G$ a cocompact hyperbolic Coxeter group and denote by $S$ the
(natural) set of generating reflections. For each generator $\,s\in S\,$, one
has $\,s^ 2=1\,$, while two distinct elements $\,s,s'\in S\,$ satisfy either no relation if the corresponding mirrors admit a common
perpendicular or provide the relation 
$\,(ss')^ {m}=1\,$
for an integer $\,m=m(s,s')>1\,$ if the mirrors intersect. 
The images of the mirrors decompose $\Bbb H^ n$
into connected components each of whose closures gives rise to a compact convex fundamental
polytope $\,P\subset \Bbb H^ n\,$ for $G$ with dihedral angles of type $\,\pi/p\,$ where
$p\ge2\,$ is an integer. Hence, $P$ is a simple polytope so that each $k$-face is contained in exactly $n-k$ facets. We call $P$ a Coxeter polytope and use the standard notation by means of  the associated
Coxeter graph simultaneously for $G$ 
and $P$ (cf. [D, Chapter 3] and [V, Chapter 5]). In particular, two nodes in the Coxeter graph $\Gamma$ of $G$ corresponding to 
mirrors intersecting under the angle of $\pi/3$ (respectively $\pi/p$) are connected by a simple 
edge (respectively by an edge with label $p$). If two mirrors are perpendicular (or admit a 
common perpendicular), their nodes are not joined at all (are joined by
a dotted line).
\smallskip
In the focus of this work is the growth series of $G$ defined by
$$f_S(x)=\sum\limits_{w\in G}x^ {l_S(w)}=1+\vert S\vert\,x+\ldots=1+\sum\limits_{i\ge 1}\,a_ix^ i\quad,$$
where $l_S(w)$ denotes the (minimal) word length of $w$ with respect to $S$, and where $a_i$
is the number of words $w$ with $l_S(w)=i$.
We investigate its explicit properties and the growth rate $\tau$ as given by the
inverse of the radius of convergence of $f_S(x)$.

In this context, the following classical facts are of fundamental importance. By a result of Steinberg [St],
$f_S(x)$ is the power series of a rational function.
For a cocompact hyperbolic Coxeter group, a result of J. Milnor [Mi] implies that $\,\tau>1\,$ 
coincides with the biggest (real) pole of $f_S(x)$ (see also [D, \S 17.1, p. 322]). Furthermore, in the same case,
the rational function $f_S(x)$ is reciprocal (resp. anti-reciprocal) for $n$ even (resp. $n$ odd)
(cf. [ChD, Corollary, p. 376] and, for $G$ having only finite Coxeter subgroups, [Se]). More precisely,
$$f_S(x^ {-1})=\cases{\,\,\,\, f_S(x)&\quad for\quad $n\equiv 0\,\,(2)\,\,,$\cr
-f_S(x)&\quad for\quad $n\equiv 1\,\,(2)\,\,.$\cr}\eqno(1.1)$$
Very useful is Steinberg's formula [St] 
$${1\over f_S(x^ {-1})}=\sum\limits_{{G_T<G\atop\scriptscriptstyle{{finite}}}}\,{(-1)^ {\vert T\vert}\over f_T(x)}\quad,\eqno(1.2) $$
allowing to express $f_S(x^ {-1})$ in terms of the growth series $f_T(x)$ of the finite Coxeter subgroups $G_T\,,\,T\subset S\,,$ of
$G$ where $\,G_\varnothing=\{1\}\,$. Recall that any subset $\,T\subset S\,$ generates a Coxeter group $G_T$ which may be finite
or infinite, reducible or irreducible. A finite Coxeter subgroup $\,G_T<G\,$ arises as stabiliser of a certain face
of $P$ and has a growth function $\,f_T(x)\,$ which, by a result of L. Solomon [So],  is a polynomial given by a product
$$f_T(x)=\prod\limits_{i=1}^ {t}\,[m_i+1]\quad.\eqno(1.3)$$
Here we use the standard notations $\,[k]:=1+x+\cdots+ x^{k-1}\,$, $\,[k,l]=[k]\cdot[l]\,$ and so on, and denote by $\,m_1=1,m_2,\ldots,m_{t}\,$
the exponents of the Coxeter group $G_T$ (cf. Table 1; for references, see [CoM, \S 9.7] or [D, Chapter 17], for example). In particular,
a maximal finite Coxeter subgroup $G_T$ of $G$ acting on $\Bbb H^ n$ is of rank $\,\vert T\vert=n\,$ and
stabilises a vertex of $P$ whose vertex neighborhood is a cone over a spherical $(n-1)$-simplex $P_v$
due to the simplicity of $P$.

\bigskip
\bigskip
\centerline{
\vbox{\offinterlineskip
\hrule
\halign{&\vrule#&
\strut\quad#\quad\cr
height6pt&\omit&&\omit&&\omit&\cr
&\hfil Graph\hfil&&Exponents\hfil&&\hfil Growth series $f_S(x)$\hfil&\cr
height4pt&\omit&&\omit&&\omit&\cr
\noalign{\hrule}
height6pt&\omit&&\omit&&\omit&\cr
&$A_n$\hfil&&${1,2,\ldots,n-1,n}$&&$[2,3,\ldots,n,n+1]$\hfil&\cr
height3pt&\omit&&\omit&&\omit&\cr
&$B_n$&&$1,3,\ldots,2n-3,2n-1$\hfil&&$[2,4,\ldots,2n-2,2n]$\hfil&\cr
height3pt&\omit&&\omit&&\omit&\cr
&$D_n$&&${1,3,\ldots,2n-5,2n-3,n-1}$&&$[2,4,\ldots,2n-2]\cdot[n]$\hfil&\cr
height3pt&\omit&&\omit&&\omit&\cr
&$G_2^{(m)}$\hfil&&$1,m-1$\hfil&&$[2,m]$\hfil&\cr
height3pt&\omit&&\omit&&\omit&\cr
&$F_4$\hfil&&$1,5,7,11$\hfil&&$[2,6,8,12]$\hfil&\cr
height3pt&\omit&&\omit&&\omit&\cr
&$E_6$\hfil&&$1,4,5,7,8,11$\hfil&&$[2,5,6,8,9,12]$\hfil&\cr
height3pt&\omit&&\omit&&\omit&\cr
&$E_7$\hfil&&$1,5,7,9,11,13,17$\hfil&&$[2,6,8,10,12,14,18]$\hfil&\cr
height3pt&\omit&&\omit&&\omit&\cr
&$E_8$\hfil&&$1,7,11,13,17,19,23,29$\hfil&&$[2,8,12,14,18,20,24,30]$\hfil&\cr
height3pt&\omit&&\omit&&\omit&\cr
&$H_3$\hfil&&$1,5,9$\hfil&&$[2,6,10]$\hfil&\cr
height3pt&\omit&&\omit&&\omit&\cr
&$H_4$\hfil&&$1,11,19,29$\hfil&&$[2,12,20,30]$\hfil&\cr
height6pt&\omit&&\omit&&\omit&\cr}
\hrule}
}
\centerline{Table 1. \it Exponents and growth polynomials of irreducible spherical Coxeter groups}
\bigskip

Finally, the growth series $f_S(x)$ of a Coxeter group acting cocompactly on $\Bbb H^n$ 
is related to the Euler characteristic of $G$ and the volume of $P$, and therefore of $\,\Bbb H^n/G\,$, as follows (see
[He]).
$${1\over f_S(1)}=\chi(G)=\cases{\displaystyle{(-1)^{n\over2}\,2\,\hbox{vol}_n(P)\over\hbox{vol}_n(\Bbb S^n)}&,\quad if $n$ is even\quad,\cr
0&,\quad if $n$ is odd\quad.\cr}\eqno(1.4)$$


Of special interest is the arithmetic nature of the growth rate $\tau>1$. By results of [Ca], [Pa] (see also [Hi]), it is known that the growth rate $\tau$  of a Coxeter group $G$ acting cocompactly on $\Bbb H^n$ is a {\it Salem number} if $\,n=2,3\,$. 
That is, $\tau>1$ is a real algebraic integer all of whose conjugates have absolute value not greater than 1, and at least one has absolute value equal to 1. It follows that the minimal polynomial of $\tau$ is palindromic with roots coming in inversive pairs.
For $n\ge4$, the growth rate $\tau$ is not a Salem number anymore. This was first observed by Cannon [Ca] who considered
$n=4$ and Coxeter groups $G$ with 5 generators. Based on substantial experimental data, we make the following claim
concerning the {\it positive} poles (appearing in inversive pairs) of the growth function of $G$
(cf. [Pe, \S 5).

\def\Vor{\parindent=17pt\par\hang\textindent}
\def\ref#1{ [#1]}

\bigskip
\medskip
{\bf Conjecture.}\quad {\sl Let $G$ be a Coxeter group acting cocompactly on $\Bbb{H}^n$ with natural
generating set $S$ and growth series $f_S(x)$. Then, }

\Vor{\ad{(a)}} {\it for $n$ even, $f_S(x)$ has precisely ${n\over 2}$ poles $0 < x_1 < \ldots < x_{{n\over2}} < 1$
in the open unit interval $]0,1[$};
\smallskip
\Vor{\ad{(b)}} {\it for $n$ odd, $f_S(x)$ has precisely the pole $1$ and ${n-1\over2}$ poles $0 < x_1 < \ldots < x_{{n-1\over2}} < 1$
in the interval $]0,1]$}.
\parindent=0pt

{\it In both cases, the poles are simple, and the non-real poles of $f_S(x)$ are contained in the annulus of radii $x_{\star}$ and  $x_{\star}^{-1}$ for some $\star \in \{1, \ldots, \left[{n\over2} \right]\}$.} 

\bigskip


In order to study such arithmetical properties of the growth series of an arbitrary cocompact Coxeter 
group $G$ acting with generating set $S$ in $\Bbb{H}^n$, we need to control the denominator polynomial of
its growth series $f_S(x)$ and assoicate to it a certain complete 
form,
$$
f_S(x) = {P(x)\over Q(x)} = {\prod\limits_{i=1}^{r} [n_i]\over
\sum\limits_{i=0}^{d} b_{i}x^{i}}\quad, 
$$
where $P(x)\,,\,Q(x)\in\Bbb Z[x]\,$ are of equal degree, and where $[k] $ is as in (1.3). The integers 
$r$ and $n_1, \ldots, n_r \geq 2$ are related to the finite Coxeter subgroups of $G$
and their exponents. Inspired by an idea of Chapovalov, Leites and Stekolshchik 
[CLS], we are able to derive a recursion formula
for the coefficients $b_i$
of the denominator polynomial $Q$. In the recursion appear beside $\vert S\vert$ and $r$ certain counting functions such as
$\,N_k=\hbox{card}\,\{\,n_l>k\,\mid\,1\le l\le r\}\,$ for $2\le k\le i\,$ related to the numerator $P$ of $f_S(x)$ (cf. (2.9) and (2.13)).
In this way, for a given group $G$, 
we dispose of an algorithm to determine the poles and the coefficients $a_i$ in the growth series
$$f_S(x) = 1+\sum\limits_{i\ge 1}\,a_ix^ i$$
in a completely explicit 
manner and to control the growth of words in the Cayley diagram of $G$ with respect to the
word metric induced by $S$. Notice that 
the cardinalities $a_i$ are usually very difficult to determine since they 
depend on the number and the relations between the generators. Hence, it is not surprising
that all our formulas depend heavily on the combinatorics of the subgroup structure of $G$ as well.

Nevertheless, there are various applications of our recursion formulas. Firstly, 
we apply the 
recursion to the elementary family of compact right-angled hyperbolic Coxeter polytopes 
and confirm our Conjecture as formulated above (see Proposition 3.2). 
For such a polytope $P$ in $\Bbb{H}^4$, having $f_0$ vertices 
and $f_3$ facets, the associated growth series is given by
$$
{(1+x)^4\over 1 + (4-f_3)x + (f_0 - 2f_3 + 6)x^2 + (4-f_3)x^3 + x^4}\quad,
$$ 
has precisely 2 inversive pairs of positive simple poles, and is, by (1.4),  of covolume equal to
$$\hbox{vol}_4(P)={f_0-4\,f_3+16\over12}\,\pi^2\quad.\eqno(1.5)$$

This result is simpler in its appearance and less specific with its consequences than a construction presented in [Z].
Therein, the growth rates of the infinite sequence
of cocompact Coxeter groups acting on $\Bbb H^ 4$ are determined which
are constructed as $m$-garlands based on
the doubly truncated Coxeter orthoschemes $[5,3,5,3]$ and $[4,3,5,3]$. The denominator of each of these growth functions is a palindromic polynomial
of degree 18 with exactly two pairs of real (simple) roots $\,x_m^{-1}<y_m^{-1}<1<y_m<x_m\,$ while all the other conjugates lie on the unit circle.
\smallskip
At the end, we shall apply our results in order to confirm our Conjecture about the growth behavior of cocompact
Coxeter groups acting with at most 6 generating reflections on $\Bbb H^4$ (see Theorem 4.1).  We shall discuss these aspects
by outlining proofs, only (cf. [Pe]).
\bigskip{\sl Acknowledgement.}\quad The authors would like to thank Alexandr Kolpakov for sharing a nice idea 
when proving Proposition 3.2 (b).
\bigskip\bigskip
\vbox{
{\bfmittel 2. Recursion formulas for growth coefficients}
\bigskip
{\bf 2.1. The complete form.}\quad
Let $G$ be a Coxeter group acting cocompactly on hyperbolic space $\Bbb H^ n$. Denote by $S$
its natural set of generating reflections, and consider the growth series of $G$,
$$f_S(x)=1+\vert S\vert\,x+\ldots=1+\sum\limits_{i\ge 1}\,a_ix^ i\quad,\eqno(2.1)$$
}
which is an (anti-)reciprocal rational function for $n$ even (odd) according to (1.1). It can be written 
as a quotient $\,f_S(x)={p(x)\over q(x)}\,$ of relatively prime polynomials $\,p,q\in\Bbb Z[x]\,$. 
By (1.2) and (1.3), the polynomials $\,p,q\,$
are of equal degree over the integers.
On the other hand, consider the denominator of the sum in Steinberg's formula (1.2)
$$
\sum\limits_{T \in {\cal F}} {(-1)^{|T|}\over f_{T}(x)}, 
$$
where ${\cal F} = \{ T \subset S\,\vert\, G_T \ \hbox{is finite}\}$. The least common multiple
$$\hbox{Virg} (S) := \hbox{LCM} \{ f_{T}(x)\,\vert\,T \in {\cal F} \}
$$
is called the {\it virgin form} of the numerator of $(-1)^n\,f_{S}(x)$, and $(-1)^n\,f_{S}(x)$ can be expressed as a rational function with numerator equal to $\hbox{Virg} (S)$ (see [CLS, Corollary 5.2.2a.]). 
Although each constituent $\,f_{T}(x)=\prod\limits_{i=1}^ {t}\,[m_i+1]\,,\,T \in {\cal F},\,$ is a product of polynomials
of type $[k]$ according to (1.3), certain factorisation properties of $[k]$ hinder $\hbox{Virg} (S)$
to be a product of $[k]$'s, only (cf. Example 1). More precisely, there is the factorisation (cf [Pr, \S 3.3])
$$[k]=\prod\limits_{d\mid k\atop d>1}\Phi_d(x)\quad, $$
where $\Phi_d(x)$ denotes the $d$-th cyclotomic polynomial of degree equal to Euler's function $\varphi(d)$. The
polynomial $\Phi_d(x)$ is irreducible in $\Bbb Z[x]$ and, for $\,d>2\,$, of even degree. If $p$ is prime and
$\,d=pm\,$, it
satisfies the  property 
$$\Phi_{pm}(x)=\cases{\Phi_m(x^p)&if $\,p\mid m\,,$\cr
\displaystyle{\Phi_m(x^p)\over\Phi_m(x)}&else~.\cr}$$
Since, for later purposes, we are interested in having uniformly tractable numerators for $f_S(x)$, we modify $\hbox{Virg}(S)$ in the following way.
Denote by $\,\hbox{Ext}(S)\in\Bbb Z[x]\,$ the monic polynomial arising as the unique common multiple
of all $\,f_{T}(x)\,,\,T \in {\cal F},\,$ such that
$$\hbox{Ext}(S)=\prod\limits_{i=1}^r\,[n_i]\quad,$$
where the integers $\,r,n_1,\ldots,n_r\ge2\,$ with $n_i=m_i+1$
are minimal. Since $\,\hbox{Ext} (S)=\hbox{Virg} (S)\cdot R(x)\,$
for some polynomial $\,R(x)\in\Bbb Z[x]\,$, $\,\hbox{Ext}(S)\,$ is called the {\it extended form} of $\hbox{Virg} (S)$. Denote by 
$$P(x):= \hbox{Ext} (S)\quad\hbox{and}\quad Q(x):=(-1)^n\,q(x)\cdot R(x)$$
the extended form of the numerator $p(x)$ and of the denominator $q(x)$ of $f_S(x)$. Then, the growth series $f_{S}(x)$ can be written as a rational function $P(x)/Q(x)$ which is called its {\it complete form}, a notion
going back to Chapovalov, Leites and Stekolshchik (see [CLS, paragraph 5.4.2]).
Let us point out that $P(x)$ and $Q(x)$ are in general no more relatively prime. An important feature of putting a growth series into its complete form is that the numerator $P$ is simply a product of polynomials $[k]$ which is of advantage when taking iterative derivatives and evaluating at $0$. The passage to the complete form does not change the number of the real poles and their localisation in the complex plane. In fact, the extension of the denominator $q(x)$ arises by multiplying it with cyclotomic polynomials of degree bigger than 1.
The next example illustrates the above procedure.
\smallskip
{\bf Example 1.\quad}
Consider the cocompact hyperbolic simplex group $G_L$ acting on $\Bbb H^4$, with set $S$ of 5 reflections related by the graph
$$\Gamma_L\,:\,\gr{5}\gd\gd\gr{4}\bullet\quad,$$
and with growth series $f_{S}(x) = p(x)/q(x)$. By means of (1.2) and the list of exponents of the subgroups involved (see Table 1), one computes 
$\,\hbox{Virg}(S) = [2, 12, 20, 30] \,\Phi_8(x)\,$. Therefore,
the complete form of $f_S(x)$ is given by the quotient of 
$\,P(x) = \hbox{Ext}(S) = [2, 8, 12, 20, 30] \,$ divided by $\,Q(x) =[4]\,q(x)\,$, where we used the decomposition 
$[8]=[4]\,\Phi_8(x)\,$.
\bigskip
{\bf 2.2. The additive nature of Ext($S$).}\quad
Write $\,f_S(x)=P(x)/Q(x)\,$ with
$$P(x)=\prod\limits_{i=1}^{r}[n_i]\quad\hbox{and}\quad Q(x)=\sum\limits_{i=0}^d\,b_ix^ i\in\Bbb Z[x]\eqno(2.2)$$
according to \S 2.1.
Since $f_S$ is (anti-)reciprocal for $n$ even (odd), and since each factor $[k]$ of $P$ (and therefore $P$ itself) is a palindromic polynomial,
satisfying the property $\,F(x)=x^{\hbox{deg}F}\,F(x^{-1})\,$, the numerator $Q$ is an (anti-)palindromic polynomial for $n$ even (odd). This means that, for $n$ odd, $Q$ satisfies the property $\,Q(x)=-x^{\hbox{deg}Q}\,Q(x^{-1})\,$.
Our aim is to derive a recursion formula for the coefficients $b_i$ of $Q(x)$. Inspired by [CLS],
we will
differentiate iteratively $\,\pm 1/f_S(x)=1/f_S(x^{-1})\,$ by means of Steinberg's formula (1.2)
and compare it - after evaluation at $x=0$ - with the corresponding expression for $\,\pm Q(x)/P(x)\,$.
Since $\,f_S(0)=1\,$ and $\,P(0)=1\,$, one has
$\,b_0=\pm b_d=1\,$. Furthermore, $\,Q^{(l)}(0)=l!\,b_l\,$. One also observes that $\,P'(0)=r\,$. 
However, by (2.2), $P(x)$ is a product of factors of type $\,[k]=1+x+\cdots+x^ {k-1}\,$ so that higher derivatives
of it become complicated expressions. The following lemma about the
{\it additive} character of $P$ is therefore very useful.
\bigskip
{\bf Lemma 2.1.}\qquad {\it Let $\,r\ge 1\,$ and $\,n_1,\ldots,n_r\ge2\,$ be integers. Then,}
$$\eqalignno{
(x-1)^ {r-1}\,&\prod\limits_{i=1}^ r\,[n_i]=[n_1+\cdots+n_r]-\sum\limits_{1\le i\le r}\,[n_1+\cdots+\widehat{n_i}+\cdots+n_r]+&(2.3)\cr
&+\sum\limits_{1\le i<j\le r}\,[n_1+\cdots+\widehat{n_i}+\cdots+\widehat{n_j}+\cdots+n_r]-\cdots+(-1)^ {r-1}\,\sum\limits_{i=1}^ r\,[n_i]\quad.\cr}$$
\bigskip
{\it Proof.}\quad We proceed by induction. Since
$$[k]=1+x+\cdots+x^ {k-1}={x^ k-1\over x-1}\quad,$$
one immediately deduces that
$$\eqalignno{
[n_1]\,[n_2]=&{x^ {n_1}-1\over x-1}\cdot{x^ {n_2}-1\over x-1}={1\over(x-1)^ 2}\,\big\{\,x^ {n_1+n_2}-x^ {n_1}-x^ {n_2}+1\,\big\}\cr
=&{1\over(x-1)^ 2}\,\big\{\,(x^ {n_1+n_2}-1)-(x^ {n_1}-1)-(x^ {n_2}-1)\,\big\}&(2.4)\cr
=&{1\over x-1}\,\big\{\,[n_1+n_2]-[n_1]-[n_2]\,\big\}\quad.\cr}$$
By means of the induction hypothesis and by using (2.4), an easy rearrangement of the terms suffices to finish the proof.
\hfill$\square$ 
\bigskip
{\bf Remark 1.}\quad It is convenient to write equation (2.3) in a more efficient way by introducing the following notation.
Let $\,X=\{\,x_1,\ldots,x_k\,\}\subseteq\{\,1,\ldots,r\,\}\,$ be a non-empty index subset, with $\,x_i< x_k\,$ if $i<k$, and write
$$n_X:=n_{x_1}+\cdots+n_{x_k}\quad.$$
Then, (2.3) can be expressed in the form
$$(x-1)^ {r-1}\,\prod\limits_{i=1}^ r\,[n_i]=\sum\limits_{\emptyset\neq X\subseteq\{\,1,\ldots,r\,\}}(-1)^{r-\vert X\vert}\,[n_X]\quad.\eqno(2.3)'$$
\smallskip
Now, by differentiating $l$-times a term $[n]$ and evaluating it at $\,x=0\,$, one obtains
$$[n]^ {(l)}(0)=l!\,\epsilon_{l}(n)\quad\hbox{where}\quad\epsilon_{l}(n):=\cases{1\,\,&if\quad$l<n\,\,$,\cr
0\,\,&if\quad$l\ge n\,\,$.\cr}\eqno(2.5)$$
The $l$-th derivative of the factor $\,1/(x-1)^ {r-1}\,$ at $\,x=0\,$ yields, for $\,l\ge1\,$, 
$$\bigg({1\over(x-1)^ {r-1}}\bigg)^{(l)}(0)=(-1)^ {r-1}\,\prod\limits_{i=0}^ {l-1}(r+i-1)\quad.\eqno(2.6)$$
\bigskip
{\bf Corollary 2.2.}\qquad {\it Let $\,n_1,\ldots,n_r\ge2\,$ be integers. Then, for $l\ge1$,}
$$\eqalignno{
\bigg(\prod\limits_{i=1}^ r\,[n_i]&\bigg)^{(l)}(0)=l!\,\sum\limits_{\emptyset\neq X\subseteq\{\,1,\ldots,r\,\}}(-1)^{\vert X\vert+1}\,\epsilon_l(n_X)+\cr
&+\sum\limits_{j=0}^{l-1}\,\bigg\{\,{l!\over (l-j)!}\,\,\prod\limits_{k=1}^{l-j}\,(r-2+k)\cdot \sum\limits_{\emptyset\neq X\subseteq\{\,1,\ldots,r\,\}}(-1)^{\vert X\vert+1}\,\epsilon_j(n_X)\,\bigg\}\,\,.&(2.7)\cr}$$
\bigskip
{\it Proof.}\quad By (2.3)', we can write
$\,\prod\limits_{i=1}^ r\,[n_i]=:u(x)\cdot v(x)\,$ with
$$u(x)={1\over(x-1)^ {r-1}}\quad\hbox{and}\quad v(x)=\sum\limits_{\emptyset\neq X\subseteq\{\,1,\ldots,r\,\}}(-1)^{r-\vert X\vert}\,[n_X]\quad.$$
Then, 
$$\bigg(\prod\limits_{i=1}^ r\,[n_i]\bigg)^{(l)}(0)=(-1)^ {r-1}\,v^{(l)}(0)+\sum\limits_{j=0}^{l-1}\,{l\choose j}\,u^ {(l-j)}(0)\,v^ {(j)}(0)\quad,$$
where, by (2.5) and (2.6), for $\,l-j\ge1\,$,
$$\eqalign{u^ {(l-j)}(0)&=(-1)^ {r-1}\,\prod\limits_{k=1}^ {l-j}\,(r-2+k)\quad,\cr
v^ {(j)}(0)&=j!\,\sum\limits_{\emptyset\neq X\subseteq\{\,1,\ldots,r\,\}}(-1)^{r-\vert X\vert}\,\epsilon_j(n_X)
\quad,\cr}$$
which allows to conclude.\hfill$\square$ 

\smallskip
{\bf Remark 2.\quad} For practical purposes, the following recursive version of Corollary 2.2 is useful (cf. [Pe, Proposition 4.5]). Denote by
$\,g_r(x):=\prod\limits_{i=1}^ r\,[n_i]\,$. Then,
$${\eqalign{g^{(l)}_r(0)=(-1)^{r+1} \,l!\,&\sum\limits_{X\subsetneq\{\,1,\ldots,r\,\}}(-1)^{\vert X\vert}\,\hbox{card}\{\,(r-\vert X\vert)\hbox{-tuples }Y\,\vert\,n_Y>r\,\}+\cr
+&\sum\limits_{j=1}^l\,{l\choose j}(-1)^{j+1}\,\prod\limits_{k=1}^j\,(r-k)\,g^{(l-j)}_r(0)\quad.\cr}}\eqno(2.8)$$

As an application of Corollary 2.2, we describe the cases $l=1,2,3$ explicitly. To this end,
consider the numbers
$$N_k=N_k(G):=\hbox{card}\,\{n_i>k\,\vert\,1\le i\le r\,\}\quad,\eqno(2.9)$$
for $\,k\in\Bbb N\,$, which satisfy $\,N_0=N_1=r\,$.
\bigskip
{\bf Corollary 2.3.}\quad {\it Let $\,n_1,\ldots,n_r\ge2\,$ be integers, and let $\,g_r(x)=\prod\limits_{i=1}^ r\,[n_i]\,$. Then, }
$$\eqalign{g'_r(0)&=r\cr
g''_r(0)&=r(r-1)+2\,N_2\cr
g^{(3)}_r(0)&=r(r-1)(r-2)+6\,(r-1)\,N_2+6\,N_3\quad.\cr}\eqno(2.10)$$
\bigskip
{\it Proof.}\quad By taking once the derivative of the product $g_r(x)$ and evaluating it at $x=0$
yields the claim, since 
$$\sum_{\emptyset\neq X\subseteq\{\,1,\ldots,r\,\}}(-1)^{\vert X\vert}=\sum\limits_{k=1}^r\,(-1)^ k\,{r\choose k}=-1\quad.\eqno(2.11)$$
Consider the second derivative $g''_r(x)$. By means of (2.7), and since
$\,\epsilon_0(n_X)=\epsilon_1(n_X)=1\,$, we obtain
$$g''_r(0)=\big\{r(r-1)+2\,(r-1)\big\}\cdot\sum_{\emptyset\neq X\subseteq\{\,1,\ldots,r\,\}}(-1)^{\vert X\vert+1}+2\,N_2+2\cdot \sum_{X\subseteq\{\,1,\ldots,r\,\}\atop \vert X\vert\ge2}(-1)^{\vert X\vert+1}\,\,.$$
By (2.11), the last term can be transformed in order to yield the desired equality. As for $g^{(3)}_r(0)$, a similar consideration based on Remark 2
allows to conclude.
\hfill$\square$
\bigskip
{\bf Remark 3.\quad}We will apply Corollary 2.2 later in the following inductive way. The $l$-th
derivative of the inverse function $\,h_r(x)=1/g_r(x)\,$ evaluated at $\,x=0\,$ can be expressed in terms of the lower order derivatives of $h_r(x)$ and $g_r(x)$
at $0$ as follows.
$$\bigg({1\over g_r(x)}\bigg)^{(l)}(0)=-\sum\limits_{j=1}^l\,{l\choose j}\,g_r^{(j)}(0)\,\bigg({1\over g_r(x)}\bigg)^{(l-j)}(0)\quad.\eqno(2.12)$$
This formula is a consequence of Leibniz' rule applied to $\,g_r(x)\cdot{1\over g_r(x)}\,$ and $\,g_r(0)=1\,$.
\bigskip
{\bf 2.3. The recursion formula.}\quad
Let us return to a cocompact hyperbolic Coxeter group $G$ with set of generating reflections $S$ and 
growth series in complete form $\,f_S(x)=P(x)/Q(x)\,$, that is, $\,P(x)=\prod\limits_{i=1}^{r}[n_i]\,$
and $Q(x)=1+\sum\limits_{i=1}^db_ix^i\,$ (see (2.2)). By Steinberg's formula (1.2),
$$
{1\over f_S(x)}=(-1)^n\sum\limits_{T \in {\cal F}} {(-1)^{|T|}\over f_{T}(x)}\quad,\eqno(1.2)'
$$
where ${\cal F} = \{ T \subset S\,\vert\, G_T<G \ \hbox{is finite}\}$, as usually.
Each finite Coxeter subgroup $G_T$ of $G$ has a growth polynomial of the form
$$f_T(x)=\prod\limits_{i=1}^{\vert T\vert}\,[1+m_i]=:\prod\limits_{i=1}^{\vert T\vert}\,[c_i]\quad,\eqno(2.13)$$
where the exponents $\,m_i=m_i(T)\,$ depend on $G_T$ as indicated in Table 1. 
Let 
$$C_k=C_k(T):=\hbox{card}\,\{c_i>k\,\vert\,1\le i\le \vert T\vert\,\}\quad,$$
and consider the set
$${\cal F}': = \{ T \subset S\,\vert\, \vert T\vert\ge2 \hbox{ and }G_T\hbox{ is finite}\}\quad.\eqno(2.14)$$
We are now ready to present formulas for the coefficients $b_1,b_2,b_3$ of $Q$. Observe that the coefficient $b_1$ 
has first been described in [CLS, Theorem 5.4.3], but by a different method. In the proof of [CLS], there is furthermore
a little flaw concerning the (non-)reciprocity of $f_S(x)$ when deriving and evaluating its inverse at $x=0$.
\bigskip
{\bf Proposition 2.4.}\quad{\it Let $G$ be a Coxeter group, with set $S$ of generating reflections, which acts
cocompactly on $\Bbb H^n$. Denote by $\,f_S(x)=P(x)/Q(x)\,$ its growth series in complete form with
$\,P(x)=\prod\limits_{i=1}^{r}[n_i]\,$
and $Q(x)=1+\sum\limits_{i=1}^db_ix^i\,$.
Then,}
$$\eqalignno{b_1=&\,r-\,\vert S\vert\quad,&(2.15)\cr
2\,b_2=&(-1)^{n+1} \, 2 \, \vert S\vert + (-1)^{n}\left( \sum\limits_{T \in {\cal F'}} (-1)^{|T|} |T| \,(|T| + 1) \right)+\cr
&+(-1)^{n+1} \,2\,\left( \sum\limits_{T \in {\cal F'}} (-1)^{|T|} \,C_2 \right)-\ r (r+1) + 2\, N_2 + 2r \,b_1\quad,&(2.16)\cr
6\,b_3=&(-1)^n \, 6 \, \vert S\vert + (-1)^{n+1} \cdot \left( \sum\limits_{T \in {\cal F'}} (-1)^{|T|} \, |T| \, (|T| + 1) \, (|T| + 2) \right)+\cr
&+\ (-1)^{n} \, 6 \, \left( \sum\limits_{T \in {\cal F'}} (-1)^{|T|} \, \left( -C_3 + (|T| + 1)\, C_2  \right) \right) +\cr
&+\, r\,(r+1) \, (r+2)+\,&(2.17)\cr
&+\, 6 \, N_3 - 6 \, (r+1) \, N_2\,+\cr
&+\, 3 \, (2N_2 - r(r+1)) \, b_1 + 6 r \, b_2\quad.\cr}$$
\bigskip
{\it Proof.}\quad 
In order to determine $b_1$, recall that $f_S(x)$ is given by (2.1) in the form
$$f_S(x) = \sum\limits_{i \geq 0} a_{i}x^{i} = 1 + \vert S\vert \, x + \sum\limits_{i \geq 2} a_{i}x^{i}\quad,$$
where $a_i > 0$, $i \geq 2$, are certain cardinalities. For example, $a_1 = \vert S\vert$. Since $\sum\limits_{i \geq 0} a_{i}x^{i} $$= {P(x)/Q(x)}\,$, it follows that
$$\left( \sum\limits_{i=0}^{d}b_{i}x^{i} \right) \, \left( 1 +  \vert S\vert\,  x + a_{2}x^{2} + \ldots \right ) = \prod\limits_{i=1}^{r}[n_i] \quad.\eqno(2.18)$$
A comparison of coefficients in (2.18) leads to $\,r=b_1+\vert S\vert\,$.

As for $b_2$, one computes by means of (1.2) and (1.3) that
$$\eqalign{
\left( {1\over f_{S}} \right)'' (0) =& (-1)^{n} \,\big\{ -2 \, \vert S\vert + \sum\limits_{T \in {\cal F'}}(-1)^{|T|} \,\left( - {f_{T}''(0)\,(f_{T}(0))^2 - 2\, f_{T}(0) \,(f_{T}'(0))^2\over (f_{T}(0))^4} \right) \big \} \cr
=& (-1)^{n} \, \big\{- 2 \,\vert S\vert + \sum\limits_{T \in {\cal F'}}(-1)^{|T|} \, \left( 2\,(f_{T}'(0))^2 - f_{T}''(0) \right)\big \}\quad.\cr}$$
By Corollary 2.3, applied to $\,f_T(x)=\prod\limits_{i=1}^{\vert T\vert}[c_i]\,,$
$$\left( {1\over f_{S}} \right)'' (0) = (-1)^{n} \,\big(- 2 \, \vert S\vert + \sum\limits_{T \in {\cal F'}}(-1)^{|T|} \, \left\{ |T| \,(|T| + 1) - 2\,C_2  \right\} \big)\quad. $$
On the other hand side,
$$
\left( {1\over f_{S}} \right)'' (0) = \left( {Q\over P} \right) '' (0)\quad.$$
Since $\,Q^{(l)}(0)=l!\,b_l\,$ and $\,b_1=\,r-\,\vert S\vert\,$, we obtain by Corollary 2.3 that
$$\left( {Q\over P}\right)''(0) = r\, (r+1) - 2\,N_2 + 2 \, b_2 - 2\,r\, b_1\quad.$$
It remains then to compare the two expressions for $\,(1/f_S)''(0)\,$ in order to obtain the desired formula. In a similar way one verifies the claim for $b_3$.
\hfill$\square$
\bigskip
\bigskip
{\bf Application.}\quad The proof of the identity (2.16) for $b_2$ above can be performed for $b_1$ as well. Combined with
(2.15), it reveals then
some information about the distribution of
the {\it finite} and {\it infinite} subgroups of $G$ which is very useful. Since we are dealing here only with cocompact
groups, any infinite Coxeter subgroup of $G$ is hyperbolic, and we deduce that
$$\sum\limits_{T\subset S\atop \vert G_T\vert<\infty}\,(-1)^ {\vert T\vert}\,\vert T\vert=(-1)^{n}\,\vert S\vert\quad;\quad\sum\limits_{T\subseteq S\atop \vert G_T\vert=\infty}\,(-1)^ {\vert T\vert}\,\vert T\vert=(-1)^{n+1}\,\vert S\vert\quad.\eqno(2.19)$$
For illustration, consider the Coxeter group $G$ given by the graph
$$\Gamma\quad:\quad \bullet\cdots\gr{p}\gr{q}\gr{r}\bullet\quad,\quad p,q,r\ge 3\,,\,{1\over p}+{1\over q}>{1\over2 }\,\,,\,\,{1\over q}+{1\over r}<{1\over2 }\,\,,\eqno(2.20)$$
which acts cocompactly on $\Bbb H^ 3$ and is generated by five reflections in the facets of a certain
simplicial prism (more precisely, a {\it simply truncated orthoscheme}) of dihedral angles $\,\pi/p,\pi/q,\pi/r\,$. Each subgraph containing $\,\,\bullet\cdots\bullet\,\,$ in (2.20)  is of infinite order.
A little computation with respect to (2.20) confirms (2.19) as follows.
$$\sum\limits_{T\subseteq S\atop \vert G_T\vert=\infty}\,(-1)^ {\vert T\vert}\,\vert T\vert=2\cdot 1-3\cdot4+4\cdot5-5\cdot 1=5\quad.$$
\medskip
In general, the coefficient $\,b_k\,,\,k\ge4\,,$ of the denominator polynomial $Q$ of $f_S(x)$
can be deduced from $\,b_1,\ldots,b_{k-1}\,$ as follows.
By means of Steinberg's formula $(1.2)'$,
$$\left( {1\over f_{S}} \right)^{(k)} (0) = (-1)^{n+k+1} \, k! \, \vert S\vert + (-1)^n \, \sum\limits_{T \in {\cal F'}} (-1)^{|T|} \, \left( {1\over f_{T}} \right)^{(k)} (0)\quad,\eqno(2.21)$$
where ${\cal F'} \,$ is given by (2.14). On the other hand, the complete form of $f_S(x)=P(x)/Q(x)$ as given by (2.2) leads to
$$\left( {1\over f_{S}} \right)^{(k)} (0) = \left( {1\over P} \right)^{(k)} (0) + \sum\limits_{j = 1}^{k-1}{k\choose j} \, j! \, b_j \, \left( {1\over P} \right)^{(k-j)}(0) + k! \, b_k\quad.\eqno(2.22)$$
By comparing (2.21) with (2.22), one derives a first formula for the coefficient $b_k$ as follows.
$$\eqalignno{k! \, b_k =& (-1)^{n+k+1} \, k! \,\vert S\vert + (-1)^{n} \, P_k^{\tau} - P_{k} + B_k\quad,\quad\hbox{where} &(2.23)\cr
P_k :=& \left( {1\over P} \right)^{(k)} (0)\quad, \cr
P_k^{\tau} := &\sum\limits_{T \in {\cal F'}} (-1)^{|T|} \, \left( {1\over f_{T}} \right)^{(k)} (0) \quad,\cr
B_k := &- \sum\limits_{j = 1}^{k-1}{k\choose j} \, j! \, b_j \, \left( {1\over P} \right)^{(k-j)} (0)=- \sum\limits_{j = 1}^{k-1} {k!\over(k-j)!} \, b_j \, P_{k-j}\quad.\cr}
$$
The different terms in (2.23) can be determined as follows.
By (2.12), we obtain the recursion
$$\eqalign{
P_{k} =& - \sum\limits_{j=1}^{k} {k\choose j} \, P^{(j)}(0) \, P_{k-j} \quad,\cr
P_k^{\tau} = &\sum\limits_{T \in {\cal F'}} (-1)^{|T| + 1} \, \left( \sum\limits_{j=1}^{k} {k\choose j} \, f_{T}^{(j)}(0) \,\left( {1\over f_{T}} \right)^{(k-j)} (0) \right) \quad.\cr}\eqno(2.24)$$
Corollary 2.2 together with (2.13) yields now similar recursion identities for both parts, that is,
$$\eqalignno{
P_{k}=& - \sum\limits_{j=1}^{k} {k!\over (k-j)!} \,\bigg(\,\sum\limits_{\emptyset\neq X\subseteq\{\,1,\ldots,r\,\}}(-1)^{\vert X\vert+1}\,\epsilon_j(n_X)+\cr
&+\sum\limits_{i=0}^{j-1}\,\big\{\,{j!\over (j-i)!}\,\,\prod\limits_{l=1}^{j-i}\,(r-2+l)\cdot \sum\limits_{\emptyset\neq X\subseteq\{\,1,\ldots,r\,\}}(-1)^{\vert X\vert+1}\,\epsilon_i(n_X)\,\big\} \bigg) \,P_{k-j} \quad,\cr
P_k^{\tau} = &\sum\limits_{T \in {\cal F'}} (-1)^{|T| + 1} \, \sum\limits_{j=1}^{k} 
{k!\over (k-j)!} \,\bigg(\,\sum\limits_{\emptyset\neq Y\subseteq\{\,1,\ldots,\vert T\vert\,\}}(-1)^{\vert Y\vert+1}\,\epsilon_j(c_Y)+&(2.25)\cr
&+\sum\limits_{i=0}^{j-1}\,\big\{\,{j!\over (j-i)!}\,\,\prod\limits_{l=1}^{j-i}\,(\vert T\vert-2+l)\cdot \sum\limits_{\emptyset\neq Y\subseteq\{\,1,\ldots,\vert T\vert\,\}}(-1)^{\vert X\vert+1}\,\epsilon_i(c_Y)\,\big\}
\left( {1\over f_{T}} \right)^{(k-j)} (0)\bigg)\quad.\cr}$$

Finally, for $B_k$, we easily derive the relation 
$$B_k = - k \, P_{k-1} \, b_1 - \sum\limits_{j=2}^{k-2} \left( {k !\over(k - j) !} \, P_{k-j}\, b_j \right) + k! \,r \, b_{k-1} \quad.\eqno(2.26)$$
Then, by plugging (2.24)--(2.26) into (2.23), the following recursion concept follows where we add for
completeness the first values according to Proposition 2.4. Recall the notations
$\,{\cal F}'$, $N_k$, 
$C_k$ and $\epsilon_k(X)$ according to (2.5), (2.9), (2.14).
\bigskip
{\bf Theorem 2.5.\quad (The recursion formula)}\qquad{\it Let $G$ be a Coxeter group with set $S$ of generating reflections acting
cocompactly on $\Bbb H^n$. Denote by $\,f_S(x)=P(x)/Q(x)\,$ its growth series in complete form with $P(x)=\prod\limits_{i=1}^{r}[n_i]\,$ and
$Q(x)=1+\sum\limits_{i=1}^db_ix^i\,$.
Then, for $\,k\ge 4\,$, and with $\,P_k=({1\over P})^{(k)}(0)\,$,}
$$\eqalignno{
b_1=&\,r-\,\vert S\vert\quad,\cr
2\,b_2=&\,(-1)^{n+1} \, 2 \, \vert S\vert + (-1)^{n}\left( \sum\limits_{T \in {\cal F'}} (-1)^{|T|} |T| \,(|T| + 1) \right)+\cr
&+\,(-1)^{n+1} \,2\,\left( \sum\limits_{T \in {\cal F'}} (-1)^{|T|} \,C_2 \right)-\ r (r+1) + 2\, N_2 + 2r \,b_1\quad,\cr
6\,b_3=&(-1)^n \, 6 \, \vert S\vert + (-1)^{n+1} \cdot \left( \sum\limits_{T \in {\cal F'}} (-1)^{|T|} \, |T| \, (|T| + 1) \, (|T| + 2) \right)+\cr
&+\ (-1)^{n} \, 6 \, \left( \sum\limits_{T \in {\cal F'}} (-1)^{|T|} \, \left( -C_3 + (|T| + 1)\, C_2  \right) \right) \cr
&+\, r\,(r+1) \, (r+2)\,+\cr
&+\, 6 \, N_3 - 6 \, (r+1) \, N_2\,+&(2.27)\cr
&+\, 3 \, (2N_2 - r(r+1)) \, b_1 + 6 r \, b_2\quad,\cr
k! \, b_k =&\, (-1)^{n+k+1} \, k! \,\vert S\vert + \cr
&\,+\sum\limits_{j=1}^{k} {k!\over (k-j)!} \,\bigg(\,\sum\limits_{\emptyset\neq X\subseteq\{\,1,\ldots,r\,\}}(-1)^{\vert X\vert+1}\,\epsilon_j(n_X)+\cr
&+\,\sum\limits_{i=0}^{j-1}\,\big\{\,{j!\over (j-i)!}\,\,\prod\limits_{l=1}^{j-i}\,(r-2+l)\cdot \sum\limits_{\emptyset\neq X\subseteq\{\,1,\ldots,r\,\}}(-1)^{\vert X\vert+1}\,\epsilon_i(n_X)\,\big\} \bigg) \,P_{k-j} \,+\cr
&+\,\sum\limits_{T \in {\cal F'}} (-1)^{n+|T| + 1} \, \sum\limits_{j=1}^{k} 
{k!\over (k-j)!} \,\bigg(\,\sum\limits_{\emptyset\neq Y\subseteq\{\,1,\ldots,\vert T\vert\,\}}(-1)^{\vert Y\vert+1}\,\epsilon_j(c_Y)+\cr
&+\,\sum\limits_{i=0}^{j-1}\,\big\{\,{j!\over (j-i)!}\,\prod\limits_{l=1}^{j-i}\,(\vert T\vert-2+l)\cdot \sum\limits_{\emptyset\neq Y\subseteq\{\,1,\ldots,\vert T\vert\,\}}(-1)^{\vert X\vert+1}\,\epsilon_i(c_Y)\big\}
( {1\over f_{T}} )^{(k-j)} (0)\bigg)-\cr
&- k \, P_{k-1} \, b_1 - \sum\limits_{j=2}^{k-2} \left( {k !\over(k - j) !} \, P_{k-j}\, b_j \right) + k! \,r \, b_{k-1} \quad.\cr}$$
\medskip
It is obvious that formula (2.27) of Theorem 2.5 depends strongly on the finite Coxeter subgroups of a given group, together with
their exponents.
For a family of hyperbolic Coxeter polytopes with fixed combinatorial structure,
the algorithm of Theorem 2.5 can be implemented into a computer program by encoding the details about all
finite irreducible Coxeter groups according to Table 1.

\bigskip
\bigskip
\vbox{
{\bfmittel 3. Growth of right-angled hyperbolic Coxeter groups}
\bigskip
Consider a hyperbolic Coxeter group $G$ with presentation
$$G = \langle S=\{ s_1, \ldots, s_k \} \ |\ (s_is_j)^{m_{ij}} = 1 \rangle\quad.$$
}
Then, $G$ is called {\it right-angled} if and only if $m_{ij} \in \{ 1, 2, \infty \}$. The terminology is justified by the fact that a fundamental polyhedron $P \subset \Bbb{H}^n$ has all dihedral angles equal to $\pi/2$ (see also [PV]). Notice that each subgroup of $G$ and all $l$-faces, $2 \leq l \leq n-1$, of $P$ are right-angled. By results of Vinberg [V], there exist no cocompact right-angled Coxeter groups in $\Bbb{H}^n$ for $n \geq 5$. For $n=2$, right-angled Coxeter polygons are realisable as long as they have at least five vertices. For $n=3$, the (compact) right-angled dodecahedron is the one with the minimal number of facets (and vertices). A beautiful example in $\Bbb{H}^4$ is the compact (regular) $120$-cell of dihedral angle $\pi/2$ whose symmetry group is generated by the reflections of the Coxeter group 
$$
\gr{5}\gd\gd\gr{4}\bullet\quad.
$$
Let $G$ be a cocompact right-angled Coxeter group, with generating set $S$, and which acts on $\Bbb{H}^n$. Hence, $n \leq 4$. Consider the growth series $f_{S}(x) = P(x)/Q(x)$ of $G$ in its complete form (see \S 2). 
Each finite Coxeter subgroup $G_T$ of $G$ is right-angled and, by (1.3),  has a growth series equal to $\, [2]^{\vert T\vert}\,$. Hence,
the numerator in its virgin form of $f_S(x)$ equals $[2]^n$, since the maximal (right-angled) subgroups in $G$ are of rank $n$. We obtain
$$
f_{S}(x) = { [2]^{n}\over Q(x)}\quad\hbox{with} \quad Q(x) = 1+\sum\limits_{i=1}^{n} b_{i}x^{i}\quad.\eqno(3.1)$$
Recall that  $\,b_{n-i}=(-1)^n\,b_i\,$ for all $\,0\le i\le [n/2]\,$, since $Q(x)$ is (anti-)palindromic. Furthermore $n \leq 4$, so that
at most the coefficients $b_0,b_1,b_2$ are of pertinence in (3.1).  As a consequence of Theorem 2.5, one deduces easily the following result.
\bigskip
{\bf Corollary 3.1.\quad}{\it
Let $G$ be a right-angled hyperbolic Coxeter group, with generating set $S$, which acts cocompactly on $\Bbb{H}^n$, $n \leq 4$. Then, the coefficients $\,b_i\,,\,1\le i\le [n/2]\,$, of $Q(x)$ in (3.1) are given by}
$$\eqalign{
b_1=&\,n-\vert S\vert\cr
b_2=&\,{n\over 2}\,(n-2\,\vert S\vert -1)\,+\,{(-1)^n\over 2}\,\bigg(\,\sum\limits_{T \in {\cal F'}} (-1)^{\vert T\vert} \,\vert T\vert (\vert T\vert+ 1)-2\,\vert S\vert\bigg)\quad,\cr}\eqno(3.2)$$
{\it where $\,{\cal F}'= \{ T \subset S\,\vert\, \vert T\vert\ge2 \hbox{ and }G_T\hbox{ is finite}\}\,$.}
\bigskip
{\bf Remark 4.}\quad Consider a group as in Corollary 3.1 together with its growth series 
$\,f_S(x)=\sum\limits_{i \geq 0} a_{k}x^{k} = 1 + \vert S\vert x + \sum\limits_{k\geq 2} a_{k}x^{k}\,$ where $a_k > 0$ is the number of $S$-words of length $k$ in $G$. Formula (3.1) yields the following recursion for $a_k$ with $a_0=1$ and $a_1=\vert S\vert\,$.
$$a_k=\cases{{n\choose k} - \sum\limits_{j = 1}^{k} a_{k-j}b_{j}&\quad for $\,2 \leq k \leq n\,\,;$\cr
- \sum\limits_{j = 1}^{n} a_{k-j}b_{j}&\quad for $\,2 \leq k \leq n\,\,,$\cr }$$
where $b_0=1$, and the coeffcients $b_i\,,\,i = 1, \ldots, n$, are given by Corollary 3.1, and  $b_i= 0$ for $i > n\,$.

For example, a {\it hexagonal} right-angled Coxeter group $G_H$ acting cocompactly on $\Bbb H^2$ has a Coxeter series 
$\,f_S(x)=1 + 6 x + 24x^2+90x^3+336x^4+1254x^5+4680x^6+17466x^7+65184x^8+243270x^9+907896x^{10}+\cdots\,$.
\bigskip
In the sequel, we present a combinatorial formula for $b_2$ in Corollary 3.1. Consider an arbitrary convex $n$-polytope $P \subset \Bbb{H}^n$. Its {\it $f$-vector} $f = (f_0,f_1,\ldots,f_{n-1})\,$ 
has components $f_i$ given by 
the numbers of $i$-faces of $P$. They are related by Euler's formula according to 
$$\sum\limits_{i=0}^{n-1} (-1)^{i} \, f_i = 1 -(-1)^n\quad.\eqno(3.3)$$
\bigskip
{\bf Proposition 3.2.\quad}{\it
Let $G$ be a right-angled Coxeter group, with generating set $S$, acting cocompactly on $\Bbb{H}^4$ with fundamental polytope $P$ and $f$-vector $\,(f_0,f_1,f_2,f_3)\,$. Let $f_S(x)$ denote the growth series of $G$ in its complete form. Then, }
\medskip
\settabs10\columns
\+(a)&$\displaystyle{f_S(x) ={[2]^4 \over 1 + (4-f_3)x + (f_0 - 2\, f_3 + 6)x^2 + (4-f_3)x^3 + x^4}\quad;}$\cr
\bigskip
\+(b)&{\it $f_S(x)$ has four distinct (simple) poles given by}\cr
\medskip
\+&$x_1=$&$\displaystyle{{1\over4} \left( \alpha + \sqrt{\gamma} + \sqrt{\beta +2\alpha\sqrt{\gamma}} \right) \quad;\quad x_1^{-1} ={1\over4} \left( \alpha + \sqrt{\gamma} - \sqrt{\beta +2\alpha\sqrt{\gamma}} \right)}\,$;\cr
\+&$x_2=$&$\displaystyle{{1\over4} \left(\alpha - \sqrt{\gamma} + \sqrt{\beta -2\alpha\sqrt{\gamma}} \right) \quad;\quad x_2^{-1} ={1\over4} \left( \alpha - \sqrt{\gamma} - \sqrt{\beta -2\alpha\sqrt{\gamma}} \right)}\,$,\cr
\medskip
\+{\it where}\cr
\+&$\alpha=\,f_3-4\quad,\quad\beta=\,2 \alpha f_3 - 4 f_0\quad,\quad\gamma=f_3^2 - 4 f_0\quad$\quad.\cr
\bigskip
\+(c)&{\it $-1$ is a root of multiplicity 4 of $f_S(x)$\quad.}\cr
\bigskip
\+(d)&{\it The volume of $P$ is given by}\cr
$$\hbox{vol}_4(P)={f_0-4\,f_3+16\over12}\,\pi^2\quad.$$
\bigskip
{\it Proof.}\quad As for (a), it is sufficient by Corollary 3.1 to show that $b_2 = f_0 - 2\, f_3 + 6$. Since $P$
is simple, $2 f_0 = f_1$, and the number of finite Coxeter subgroups $G_T$ of $G$ with $\vert T \vert = l$ equals $f_{4-l}$, for $l = 1, \ldots, 4$. 
Moreover, $f_3 = \vert S\vert$, and by Euler's formula (3.3), $f_2 = f_0 + \vert S\vert$. Hence, by Corollary 2.6,
$$\sum\limits_{T \in {\cal F'}} (-1)^{\vert T\vert} \, \vert T\vert\,(\vert T\vert+ 1) = 6 f_2 - 12 f_1 + 20 f_0=6 \,\vert S\vert + 2 \, f_0=6\,f_3+ 2 \, f_0\quad.\eqno(3.4)$$
By plugging (3.4) into (3.2), we deduce that
$\,b_2 = f_0 - 2 \, f_3+ 6\,$. Property (c) follows easily from (a) since $\,f_S(-1)=\displaystyle{[2]^4(-1)\over f_0}=0\,$,
and property (d) is a direct consequence of (a) and Heckman's formula
(1.4).

As for (b), consider the denominator $Q(x) = 1 + (4 - f_3)x + (f_0 - 2f_3 + 6)x^2 + (4-f_3)x^3 + x^4$ of $f_S(x)$ in (a). The polynomial $Q(x)$ is quartic over the integers with discriminant (see [R, {\it Discriminants}])
$$\Delta = f_0 \, (16 + f_0 - 4\, f_3) \, (f_3^2 - 4 \, f_0)^2\quad.$$
Since $P$ is simple with 2-faces being at least pentagonal, $f_0 \geq 5f_3$. Suppose that $f_3^2 > 4f_0$. Then, $\Delta >0$, and  $Q(x)$ 
has only simple, real roots. It is a standard matter to determine the explicit form of these roots. In fact, 
by applying the transformation $\,x = X + 1/X\,$ to the quartic polynomial $Q(x)$, which does not change 
the discriminant, one obtains a reduced cubic $\widetilde{Q}(x)$ with explicit formulas for its roots 
(see [R, {\it Classical Formulas}]). It remains to show that $f_3^2 > 4f_0$. Suppose on the contrary 
that $f_3^2 \le 4f_0$ and consider
a facet $F^*$ of $P$ with maximal number $\,N:=f_0(F^*)=\hbox{max}\,f_0(F)\,$ of vertices among all facets $F$ of $P$, that is,
$$4f_0=\sum\limits_{F}f_0(F)\le Nf_3\quad\hbox{whence}\quad f_3\le N\quad.$$
We conclude the proof by showing that $\,N\ge f_3\ge f_0(F^*)+1=N+1\,$.
Indeed, since $P$ is simple, precisely one additional edge of $P$ emanates from each of the $\,N=f_0(F^*)$
vertices
to the outside of $F^*$.
We show that these $N$ edges give rise to $N$ different (but not necessarily disjoint) facets $F_*$
of $P$, beside $F^*$, and this by contraposition. Since all facets are convex and meet properly at 2-faces of $P$, 
the assumption of the opposite can hold only if two vertices $v_1,v_2$ belong to a common edge
of $F^*$ and if their edges leaving $F^*$ 
arrive at vertices $w_1,w_2$ which may coincide in or lie on an edge
of $F_*$. Hence, the
convex hull of the vertices $\,v_1,v_2,w_1,w_2\,$ is a right-angled triangular or quadrilateral 2-face of $P$ which is impossible.
Therefore, $\,N\ge f_3\ge f_0(F^*)+1=N+1\,$.
\hfill{$\square$}
\bigskip
{\bf Remark 5.\quad}In [D, Example 17.4.3], a result analogous to Proposition 3.2 (a) for the $3$-dimensional case is presented. More precisely, the growth series of a cocompact right-angled Coxeter group in $\Bbb{H}^3$ is given by 
$$
f_S(x) ={[2]^3\over1 - (f_0 - 3)x + (f_0 - 3)x^2 - x^3}\quad,\eqno(3.5)
$$ 
which has the three positive real poles $1$, $\tau$ , $\tau^{-1}$ where
$$\tau = {(f_0 - 4) + \sqrt{(f_0 - 4)^2 - 4}\over2}\quad.$$
In fact, Euler's formula (3.3), $\,f_0-f_1+f_2=0\,,$ together with the (vertex) simplicity $\,3\,f_0=2\,f_1\,$ and
the evident inequality $\,f_2\ge4\,$, yields that $\,f_S(x)\,$ has only real simple roots.
\medskip
{\bf Example 2.\quad}
Let $G_{120}$ be the Coxeter group generated by the $120$ reflections with respect to the facets of a right-angled (compact)  $120$-cell $P \subset \Bbb{H}^4$. The polyhedron $P$ has $f$-vector
$$
f = \left( 600, 1200, 720, 120 \right) 
$$
and is the $4$-dimensional analogue of a right-angled dodecahedron $D$. In fact, all facets of $P$ are isometric to $D$. 
The volume of $P$ equals
$\,34\pi^2/3$ by Proposition 3.2 (d). This value can also be obtained by studying the symmetry group of $P$ and by determining the covolume of its 14'000 index Coxeter simplex subgroup according to [Ke, Appendix], that is,
 $$\hbox{vol}(P)=14'400\cdot\hbox{covol}_4(\,\gr{5}\gd\gd\gr{4}\bullet\,)=14'400\cdot{17\pi^2\over21'600}={34\pi^2\over3}\quad.$$
By means of (3.1) and Proposition 3.2, the growth series of $G_{120}$ with respect to the set $S$ of the above reflections is given by 
$$
f_{S}(x) = {[2]^4\over1 - 116 x + 366x^2 - 116x^3 + x^4}\quad,\eqno(3.6)$$
implying that $f_S(x)$ possesses exactly two pairs of real poles, which are positive and simple.
\bigskip\bigskip
{\bf 4. Growth of Coxeter groups with at most 6 generators in $\Bbb H^ 4$}
\bigskip
Consider a hyperbolic cocompact Coxeter group $G$ with generating set of reflections $S$ acting in low dimensions $n\ge2$.
For $n=2$, J. Cannon and P. Wagreich [CaW] showed that the growth series $f_S(x)$ is a quotient of
relatively prime monic polynomials over the integers for which the denominator splits into exactly
one Salem polynomial and (possibly none) distinct irreducible cyclotomic polynomials.
Here, a Salem polynomial is a palindromic irreducible monic polynomial over the integers
with exactly one (inversive) pair of real roots $\,\alpha^ {-1},\alpha>1\,$ and with all other conjugates lying
on the unit circle. The root $\alpha$ is called a Salem number.
Hence,  the growth rate $\tau$ of any (of the infinitely many) planar cocompact hyperbolic Coxeter groups
is a Salem number. In [Hi], E. Hironaka showed that
the smallest growth rate which arises in this way 
equals Lehmer's number given by the root $\,\alpha_L>1\,$ of the Salem polynomial of smallest known degree
$$L(x)=1+x-x^ 3-x^ 4-x^ 5-x^ 6-x^ 7+x^ 9+x^{10}\quad.$$
For $n=3$, W. Parry [Pa] proved a result, analogous to the one in [CaW] above, for cocompact hyperbolic Coxeter groups  
providing a unified proof for $n=2$ and $n=3$ and extending Cannon's result [Ca] from the case
of the nine Coxeter tetrahedra to arbitrary compact Coxeter polyhedra. Parry's proof
is based on a special relationship between anti-reciprocal functions and Salem numbers.
However, for $n\ge4$, growth rates of cocompact hyperbolic Coxeter groups are not Salem numbers anymore
as is illustrated by the example of the compact right-angled 120-cell in $\Bbb H^4$ according to (3.6). 
\smallskip
In the following, we will describe in detail (see Theorem 4.1) the growth series of a cocompact hyperbolic Coxeter group, generated by
at most six reflections in $\Bbb H^4$, and show that its
positive poles arise always in precisely 2 inversive pairs $\,x_1^{-1}<x_2^{-1}<1<x_2<x_1\,$, and each one is of multiplicity 1.
The growth rate $\tau=x_1$ is a {\it Perron number}, that is, $\tau$ is a real algebraic integer all of whose conjugates are of strictly smaller absolute value.
The non-real poles of the growth function come in quadruplets which
do not all lie on the unit circle anymore (cf. also [CLS]).
A rigorous proof of all our observations is very technical (and partially computer-based) and necessitates a closer analysis 
of the Coxeter groups under consideration
with respect to their Coxeter subgroup structure (see [Pe] for all details). 
Indeed, here lies the qualitative 
difference to the lower dimensional
cases $\,n=2\,$ resp. $\,n=3\,$ where the maximal finite Coxeter subgroups are dihedral groups resp. spherical 
triangle groups with a very limited, manageable variety of exponents, and this independently of the number of generators.
The following exposition will document that similar growth questions in higher dimensional hyperbolic spaces become nearly
intractable.
\medskip
Let $G$ be a Coxeter group, with natural generating set $S$ such that $\vert S\vert\le6$, and which acts
cocompactly on $\Bbb H^4$.
There is a complete classification which shows that they
fall combinatorially into three (finite) families. For $\,\vert S\vert=5\,$, $G$ is a Coxeter simplex group, denoted by $G_L$. They
were discovered and classified by F. Lann\'er (cf. [V, Chapter 3, Table 3) and are nowadays called {\it Lann\' er groups}.
If $\,\vert S\vert=6\,$, then $G$ is a {\it Kaplinskaya group} $G_K$
or an {\it Esselmann group} $G_E$, which are characterised as follows. A 
fundamental polytope of $G_K$ is a product of a 1-simplex with a 3-simplex and has eight vertices, while a
fundamental polytope of $G_E$ is a product of two triangles with nine vertices. The classification 
of the Kaplinskaya groups can be found in [K], and the list of all Esselmann groups is in [E].

\vbox{
\beginpicture

\setcoordinatesystem units <1.5pt,1.5pt>
\setplotarea x from -150 to 60, y from -60 to 40

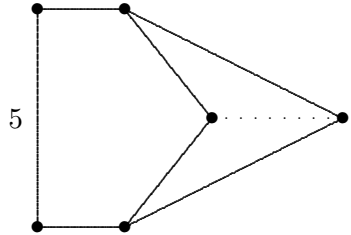
\put  {Figure 1.\quad \it The graph of the Kaplinskaya group $G_{66}$} at -10 -50
\setlinear
\plot -10 25  10 0  -10 -25  40 0  -10 25  -30 25  -30 -25  -10 -25 /

\put {$5$} at -35 0
\put {$\bullet$} at -10 25
\put {$\bullet$} at 10 0
\put {$\bullet$} at -10 -25
\put {$\bullet$} at 40 0
\put {$\bullet$} at -30 25
\put {$\bullet$} at -30 -25

\setlinear
\linethickness=.4pt
\setdots
\plot 10 0  40 0 /

\endpicture
}
\bigskip
{\bf Theorem 4.1.}\quad {\it Let $G$ be a Lann\' er group, an Esselmann group or a Kaplinskaya group, respectively,
acting with natural generating set $S$ on $\Bbb H^4$. Then,
\settabs14\columns
\smallskip
\+$(1)$&the growth series $f_S(x)$ of $G$ is a quotient of relatively prime, monic and \cr
\+&  palindromic polynomials of equal degree over the integers. \cr
\smallskip
\+$(2)$&The growth series $f_S(x)$ of $G$ possesses four distinct positive real poles appearing in \cr
\+&pairs $(x_1, x_1^{-1})$ and $(x_2, x_2^{-1})$ with $x_1 < x_2 < 1 < x_2^{-1} < x_1^{-1}\,$; these poles are simple. \cr
\medskip
\+$(3)$ & The growth rate $\tau = x_1^{-1}$ is a Perron number. \cr
\smallskip
\+$(4)$ & The non-real poles of $f_S(x)$ are contained in an annulus of radii $x_2$,  $x_2^{-1}$ around the\cr
\+& unit circle. \cr
\smallskip
\+$(5)$ & The growth series $f_S(x)$ of the Kaplinskaya group $G_{66}$ with graph $K_{66}$ (cf. Figure 1)\cr
\+& has four distinct negative and four distinct positive simple real poles; for $G \neq G_{66}$, \cr
\+& $f_S(x)$ has no negative pole.\cr
}
\bigskip
\bigskip
{\bf Remark 6.}\quad The exceptional role (5) of the Kaplinskaya group $G_{66}$, having a growth series
with 4 inversive pairs of distinct real poles, was first discovered by
T. Zehrt (cf. [Z]).
\medskip
{\it Sketch of the proof of Theorem 4.1.\quad}We will only discuss the ingredients of the proofs for (1)--(3), and this especially
for the simplex case. At the end, we indicate how the proof extends for
the families $G_E$ and $G_K$ (for more details, see [Pe]). Consider a Lann\'er group
$G_L$ with natural generator set $S$ and denote by $G_T\,,\,T\subset S\,,$ a maximal finite Coxeter subgroup of $G_L$.
Associate to $G_T$ the help-function
$$ h^{L}_{T}(x) := - {1\over x+1} + {1\over3}\,\sum\limits_{U}{1\over f_{U}(x)} - {1\over2}\,\sum\limits_{V}{1\over f_{V}(x)} + {1\over f_{T}(x)}\quad,\eqno(4.1)$$
where $U$ varies over the six $2$-element subsets and $V$ varies over the four $3$-element subsets of $T$. By means of Steinberg's formula (1.2)', 
$${1\over f_{S}(x)} = 1 + \sum\limits_{G_T \hbox{ maximal}} h_T^{L}(x)\quad.\eqno(4.2)$$
Note that there exists only a very limited number of different maximal finite Coxeter subgroups in $G_L$ as the weights of the Coxeter graph of $G_L$ are at most equal to $5$. By taking into account their reducibility properties, the following 
important auxiliary result can be shown by a case-by-case study (cf. [Pe, Lemma 3.8]).
\bigskip
{\bf Lemma 4.2.}\quad {\it 
The help-function $h_T^L(x)$ (4.1) associated to a maximal finite Coxeter subgroup $G_T$ of a Lann\'er group $G_L$ can be written as the quotient
$$
h_T^{L}(x) = -x \, {n(x)\over d(x)}\quad,\eqno(4.3)$$
where $n(x)$ and $d(x)$ are palindromic polynomials of even degrees over the integers. Moreover, $d(x)$ is cyclotomic with $\hbox{\rm deg}\,d= \hbox{\rm deg}\,n + 2$. Furthermore, $h_T^{L}(x)$ is negative for $x > 0$ and strictly decreasing on $(-1 , 0)$. }
\medskip
By plugging (4.3) into (4.2), $f_S(x)$ becomes a quotient of palindromic integer polynomials of equal (even) degree.
It is a well-known result that each palindromic integer polynomial of degree $m$ can be factored
into a product of a constant times linear (if $m$ is odd), quadratic and quartic palindromic polynomials with real coefficients
(see also [Pe, Proposition D.11]). It follows from this that $f_S(x)$ is a quotient of monic palindromic integer polynomials which are prime.
\smallskip
For the study of the real poles of $f_S(x)$, it suffices to consider the interval $\,[-1,1]\,$ as $f_S(x)$ is reciprocal.
Recall that $f_S(0)=1$ and that $\,f_S(1)>0\,$ by (1.4). Furthermore, we prove the following.
\bigskip
{\bf Lemma 4.3.}\quad {\it Let $f_S(x)$ be the growth series of a cocompact Coxeter group $G$ acting on $\Bbb H^4\,$
with natural generator set $S$ satisfying $\,\vert S\vert \le 6\,$.
Then, $\,f_S(-1)=0\,$.}
\bigskip{\it Proof.\quad} 
By Steinberg's formula (1.2), it is sufficient to show that the growth polynomial $f_T(x)$ of at least one
maximal (or, equivalently, rank 4) finite Coxeter subgroup $G_T$ of $G$ factorises according to $\,f_T(x)=[2]^4\,g_T(x)\,$ where
$\,g_T\in\Bbb Z[x]\,$ with $\,g_T(-1)\neq0\,$. Since the natural generator set $S$ of $G$ is of cardinality at most 6, $G$ contains
at least one subgroup $G_T$ of type $B_4$, $D_4$, $F_4$ or $H_4$. This is due to its combinatorial 
structure (cf. also [V], [E] and [K]). By Table 1, each of the groups $B_4$, $D_4$, $F_4$, $H_4$ has only {\it odd} exponents and therefore a growth polynomial $f_T(x)$
splitting into 4 factors of type $[2k]$. Now, observe that $\,[2k]=[2]\,\sum\limits_{i=0}^{k-1}x^{2i}\,$ so that $\,[2k](-1)=0\,$.
\hfill{$\diamond$}
\medskip
All the above observations  together with Lemma 4.2 and Lemma 4.3 allow us to conclude that $f_S(x)$ is positive and strictly increasing on
$(-1,0]\,$. Since $f_S(x)$ is reciprocal, it is non-singular on $\Bbb R_{\le0}\,$.
For the study of the behavior of $f_S(x)$ on $\,I:=[0,1]\,$, we know
that $f_S(x)$ is a rational function and has a real pole $0<x_1<1$ given by the convergence radius. This follows since  coefficients $a_i$ of the series $f_S(x)$ are positive and real (cf. \S 1; [D, \S 17.1]). In particular, $x_1$ is a real algebraic integer
whose inverse $x_1^{-1}$ is the growth rate $\tau$ of $G_L$.
Hence, $\tau$ is a Perron number.

For the proof of the remaining claims, a distinction of several cases 
and the help of a computer are needed to control the graphs of the
help-functions in the decomposition (4.2) on $I$ (see [Pe, pp. 28--47]). By doing this, it turns out that 
their sum 
$$H^L(x):=\sum\limits_{G_T \hbox{ maximal}} h_T^{L}(x)={1\over f_S(x)}-1\eqno(4.4)$$
is negative on $I$, and that $H^L(x)$ is either strictly decreasing on $I$ or possesses
exactly one negative minimum in $I$.
Since $x_1$ is a pole of $f_S(x)$, and $1/f_S(1) > 0$, it follows that $H^L(x_1) = -1$ and $H^L(1) > -1$. Therefore, $H^L(x)$ can not be strictly decreasing on $I$, but possesses exactly one negative minimum $M$. That is, there is a unique $x_M \in I$ such that $H^L(x_M) = M$. Obviously, $x_M \geq x_1$, since $x_1$ equals the radius of convergence of $f_S(x)$. Summarising, we can deduce that $f_S(x)$ possesses exactly two simple poles in $I$ if $x_M > x_1$, or it has a pole of (positive) even order in $I$ if $x_M = x_1$. 

As for the simplicity of the poles $x_1,x_2$ of $f_S(x)$, there are no criteria known to us allowing to conclude it without precise knowledge
of the denominator coefficients
(cf. [Ma], for example). 
By means of the recursion formula for these coefficients (see Theorem 2.5), the computer implementation of this algorithm helps to prove
this last claim.
\medskip
In the cases of Esselmann groups $G_E$ and Kaplinskaya groups $G_K$, our strategy is essentially the same, apart from some particularities which have to be dealt with carefully. Furthermore, the help-functions
have to be adapted to the different combinatorial features of $G_E$ and $G_K$. We finish this outline by providing their explicit shapes.

Recall that an Esselmann polytope has the combinatorial type of a direct product of two triangles and possesses therefore precisely nine vertices. The Coxeter graph $\Gamma_E$ of an Esselmann group $G_E$ contains two disjoint Lann\'er diagrams, called $L_1$ and $L_2$, each of them with three nodes. Let $G_T$ be one of the nine maximal finite Coxeter subgroups of $G_E$  where $T$ denotes its natural generating set. The help-function for $h^E_T$ is defined by
$$
h_T^{E}(x) := h_T^{L}(x) + {1\over 3 \, (1+x)} - {1\over12} \,\sum\limits_{W}  {1\over f_{W}(x)}\quad,\eqno(4.5)$$
where $h_T^L(x) $ is the function (4.1) given by
$$ h^{L}_{T}(x) = - {1\over x+1} + {1\over3}\,\sum\limits_{U}{1\over f_{U}(x)} - {1\over2}\,\sum\limits_{V}{1\over f_{V}(x)} + {1\over f_{T}(x)}\quad,$$
where $U$ is a $2$-element subset, $V$ is a $3$-element subset, and where $W$ is a subset of $T$ satisfying the following condition. The set $W$ consists of four pairs of generators $\,( s_p , s_q)\,$ such that the node in $\Gamma_E$ corresponding to $s_p$ belongs to $L_1$, while the node corresponding to $s_q\,$ belongs to $L_2$. 

A Kaplinskaya polytope has the combinatorial type of a simplicial prism and possesses therefore precisely eight vertices. The Coxeter graph $\Gamma_K$ of a Kaplinskaya group $G_K$ contains a Lann\'er diagram $L$ with four nodes which represents a tetrahedron  $P$, and two additional nodes which represent  the reflections through the top respectively the bottom of the simplicial prism $P \times [0, 1]$. Let $G_T$ be one of the eight maximal finite Coxeter subgroups of $G_K$ where $T$ denotes its natural generating set. The help-function $h^K_T$ is defined by
$$h_T^{K}(x) := h_T^{L}(x) + {1\over 4 \, (1+x)} - {1\over12}\,\sum\limits_{W}{1\over f_{W}(x)}\quad,\eqno(4.6)$$
where $h_T^L (x)$ is the function (4.1) and where $W$ is a subset of $T$ containing three pairs $\,(s_{L_1} , s_b)$, $(s_{L_2} , s_b)\,$ and $\,(s_{L_3} , s_b)\,$ such that $s_{L_j}$ belongs to $L$, for $j = 1, 2, 3$, while $s_b \not\in L$.
\hfill{$\square$}

\bigskip
{\bf Remark 7.\quad}For a cocompact hyperbolic Coxeter group, acting on $\Bbb H^4$ with a set $S$ of generating reflections such that
$\,\vert S\vert\le6\,$, the growth series $f_S(x)$ can be put into the form $R(x)/S(x)$
with monic palindromic polynomials $\,R,S\in\Bbb Z[x]\,$, $\,\hbox{deg}\,R=\hbox{deg}\,S$, and
$$ R(x)=\cases{[2,8,12,20,30]&if $G$ is a Lann\'er group\quad,\cr
[2,6,8,12,20,30]&if $G$ is an Esselmann or a Kaplinskaya group\quad.\cr}\eqno(4.7)$$

For the coefficients $\beta_k\,$ of the denominator $\,S(x)=1+\sum_{k\ge1}\beta_kx^k\,$, the recursion
of Theorem 2.5 applies as well. Observe that the numerator $R$ and the denominator $S$ are in general not prime. The result (4.7) follows easily by passing to the complete form $P/Q$ and by extending $P,Q$ simultaneously
in a suitable way according to
the Coxeter subgroup structure of all $G$ of fixed type $L$, $E$ or $K$ as described in [V], [E] and [K], and by using divisibility properties of the associated exponents.
\bigskip
{\bf Remark 8.\quad}The growth series of any hyperbolic cocompact Coxeter group $G$ acting on $\Bbb H^n\,,\,n\ge2\,,$ vanishes at $-1$ {\it if at least one of its maximal finite Coxeter subgroups has a growth polynomial with odd exponents, only}. Notice that
this condition, for $n=2$, excludes only the triangle groups $\,G=(p,q,r)\,$ with $\,p,q,r\,$ odd. For $n=3$, Parry's formulas
[Pa, (0.5), (0.6)] for $f_S(x)$ imply that $\,f_S(-1)=0\,$ without imposing any condition. However, for $\,n\ge 4\,$, the above condition is in general not true anymore.
In fact, the Tumarkin group $G_*$ acting cocompactly on $\Bbb H^6$ with Coxeter graph given in Figure 2 (cf. [T])
does not have a maximal Coxeter subgroup all of whose exponents are odd. Nevertheless, its growth series splits the factor $[2]^4$ 
so that $-1$ is a root of multiplicity $4$. As a byproduct, the computation shows also, by (1.4), that the covolume of $G_*$ is given by $2'077\pi^3/17'010'000\simeq0.003786\,$.

Finally, by analysing {\it all known} examples of cocompact hyperbolic Coxeter acting in dimensions bigger than two (and less than nine), we see that 
$-1$ is always a root of the growth series !

\vbox{
\beginpicture

\setcoordinatesystem units <1.5pt,1.5pt>
\setplotarea x from -120 to 60, y from -60 to 40

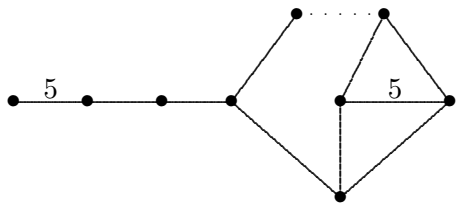
\put  {Figure 2.\quad \it A compact Tumarkin polytope with 9 facets in $\Bbb H^6$} at 10 -40
\setlinear
\plot -50 0  0 0  15 20  /
\plot 25 0  35 20  50 0  25 -22  0 0 /
\plot 50 0  25 0  25 -22 /

\put {$5$} at -41.5  3
\put {$5$} at 37.5 3
\put {$\bullet$} at -50 0
\put {$\bullet$} at -16 0
\put {$\bullet$} at -33 0
\put {$\bullet$} at 0 0
\put {$\bullet$} at 25 0
\put {$\bullet$} at 50 0
\put {$\bullet$} at 25 -22
\put {$\bullet$} at 35 20
\put {$\bullet$} at 15 20

\setlinear
\linethickness=.4pt
\setdots
\plot 15 20  35 20 /

\endpicture
}
\vskip1cm

\def\Vor{\parindent=17pt\par\hang\textindent}
\def\ref#1{ [#1]}

{\bfmittel 5. Bibliography}
\baselineskip=16pt
\bigskip
{\refklein
\Vor{\ad {[Ca]}} J. W. Cannon, The growth of the closed surface groups and the compact hyperbolic Coxeter groups,
preprint.

\smallskip
\Vor{\ad {[CaW]}} J. W. Cannon, P. Wagreich, Growth functions of surface groups, Math. Ann. {\bf 293} (1992), 239--257.

\smallskip
\Vor{\ad {[CLS]}} M. Chapovalov, D. Leites, R. Stekolshchik,
The Poincar\'e series of the hyperbolic Coxeter groups with finite volume of fundamental domains,
preprint, arXiv:0906.1596, June 2009.

\smallskip
\Vor{\ad {[ChD]}} R. Charney, M. W. Davis, Reciprocity of growth functions of Coxeter groups, Geom. Dedicata {\bf 39} (1991), 373--378.

\smallskip
\Vor{\ad {[CoM]}} H. S. M. Coxeter, W. O. J. Moser, Generators and relations for discrete groups, Springer-Verlag, 1980.  

\smallskip
\Vor{\ad {[D]}} M. W. Davis, The geometry and topology of Coxeter groups, London Mathematical Society Monographs, Princeton University Press, 2008.

\smallskip
\Vor{\ad {[E]}} F. Esselmann, The classification of compact hyperbolic Coxeter 
$d$-polytopes with $d+2$ facets, Comment. Math. Helv.  {\bf 71}  (1996), 229--242.

\smallskip
\Vor{\ad {[He]}} G. J. Heckman, The volume of hyperbolic Coxeter polytopes of even dimension, Indag. Math. (N.S.)  {\bf 6}  (1995), 189--196.

\smallskip
\Vor{\ad {[Hi]}} E. Hironaka, The Lehmer polynomial and pretzel knots, Can. Math. Soc. Bulletin {\bf 44} (2001), 440--451.

\smallskip
\Vor{\ad {[K]}} I. M. Kaplinskaja, Discrete groups generated by reflections in the faces of simplicial
prisms in Lobachevsjian spaces, Math. Notes {\bf 15} (1974), 88--91.

\smallskip
\Vor{\ad {[Ke]}} R. Kellerhals, On Schl\"afli's reduction formula, Math. Z.  {\bf 206} (1991), 193--210.

\smallskip
\Vor{\ad {[Ma]}} M. Marden, Geometry of polynomials, 2nd edition, American Mathematical Society, Providence (Rhode Island), 1966.

\smallskip
\Vor{\ad {[Mi]}} J. Milnor, A note on curvature and fundamental group, J. Differ. Geom. {\bf  2} (1968), 1--7.

\smallskip
\Vor{\ad {[Pa]}} W. Parry,
Growth series of Coxeter groups and Salem numbers,
J. Algebra {\bf 154} (1993),
406--415.

\smallskip
\Vor{\ad {[Pe]}} G. Perren, Growth of cocompact hyperbolic Coxeter groups and their rate, Ph.D. thesis, University of Fribourg, 2009.

\smallskip
\Vor{\ad {[PV]}} L. Potyagailo, E. B. Vinberg, On right-angled reflection groups in hyperbolic spaces, Comment. Math. Helv.  {\bf 80} (2005), 63--73.

\smallskip
\Vor{\ad {[Pr]}} V. V. Prasolov, Polynomials, Springer-Verlag, 2004.

\smallskip
\Vor{\ad {[R]}} J. Rotman, Galois theory, Springer-Verlag, 1998.

\smallskip
\Vor{\ad {[Se]}} J.-P. Serre,
Cohomologie des groupes discrets,
Prospects Math., Ann. Math. Stud. 70, 77--169, 1971.

\smallskip
\Vor{\ad {[So]}} L. Solomon, The orders of the finite Chevalley groups, J. Algebra {\bf 3} (1966), 376--393.

\smallskip
\Vor{\ad {[St]}} R. Steinberg, Endomorphisms of linear algebraic groups,
Mem. Amer. Math. Soc. {\bf 80} (1968).

\smallskip
\Vor{\ad {[T]}} P. Tumarkin, Compact hyperbolic Coxeter $n$-polytopes with $n+3$ facets,
Electron. J. Combin. {\bf 14} (2007), 36 pp.

\smallskip
\Vor{\ad {[V]}} E. B. Vinberg, The absence of crystallographic groups of reflections in Lobachevsky spaces of large dimensions, Trans. Moscow Math. Soc. {\bf 47} (1985), 75--112.

\smallskip
\Vor{\ad {[Z]}} T. Zehrt, C. Zehrt-Liebend\"orfer, The growth function of Coxeter garlands in $\Bbb H^{4} \,$, preprint 2009.

\bye